\documentclass[11pt]{amsart}
\usepackage{amssymb,latexsym,amsmath}
\usepackage{bbm}
\usepackage{graphicx}

\DeclareMathAlphabet{\Ma}{U}{msa}{m}{n}
\DeclareMathAlphabet{\Mb}{U}{msb}{m}{n}
\DeclareMathAlphabet{\Meuf}{U}{euf}{m}{n}
\DeclareSymbolFont{ASMa}{U}{msa}{m}{n}
\DeclareSymbolFont{ASMb}{U}{msb}{m}{n}

\def\mr #1.{\mathrm{#1\,}}
\def\mrt #1.{\mathrm{\mbox{\tiny #1\,}}}
\def\mt #1.{{\mbox{\tiny $#1$}}}
\def\ms #1.{{\mbox{\small $#1$}}}
\def\ol #1.{\overline{#1}}

\def\C{\mathbb{C}}    \def\R{\mathbb{R}}
\def\Z{\mathbb{Z}}  \def\N{\mathbb{N}}

\def\Wort#1{\mbox{\fontfamily{cmr}\selectfont\mdseries\upshape #1}}

\def\1{\mathbbm 1}
\def\EINS{\1}

\def\1{\mathbbm 1}
\def\C{\Mb{C}}
\def\ot #1.{{\got{#1}}}
\def\got#1{\Meuf{#1}}
\def\al #1.{{\mathcal{#1}}}

\theoremstyle{plain}            
 \newtheorem{theorem}{Theorem}[section]
 \newtheorem*{maintheorem*}{Main Theorem}
 \newtheorem*{maintheorem.}{Main Theorem~1'}
 
 \newtheorem{proposition}[theorem]{Proposition}
 \newtheorem{lemma}[theorem]{Lemma}

 \theoremstyle{definition}       
 \newtheorem{definition}[theorem]{Definition}

 \theoremstyle{remark}
 \newtheorem{remark}[theorem]{Remark}
 \newtheorem{example}[theorem]{Example}
 
 \newtheorem{exercise}[theorem]{Exercise}

\def\Wort#1{\mbox{\fontfamily{cmr}\selectfont\mdseries\upshape #1}}

\def\KIn#1{{\fontsize{6pt}{8pt}\selectfont\Wort{#1}}}

\newfont{\Kcal}{cmsy6 scaled 1000}
\newfont{\Kgot}{eufm6 scaled 1000}

\def\Kbegin{\begin{equation} \left. \begin{array}{rcl}}
\def\Kend{\end{array} \right\} \end{equation}}

\DeclareMathSymbol{\hsemi}{\mathord}{ASMb}{"6E}
\newcommand{\semi}[2]{\mbox{$#1\kern.1em\hsemi\kern.1em#2$}}

\def\vplatz#1{{\rule{0mm}{#1}}}

\def\LA{\left\langle\bgroup}
\def\LE{\left[\bgroup}
\def\LG{\left\{\bgroup}
\def\LR{\left(\bgroup}
\def\RA{\egroup^{\rule{0mm}{0mm}}\right\rangle}
\def\RE{\egroup^{\rule{0mm}{2mm}}\right]}
\def\RG{\egroup^{\rule{0mm}{2mm}}\right\}}
\def\RR{\egroup^{\rule{0mm}{2mm}}\right)}
\def\Ldummy{\left.\bgroup}
\def\Rdummy{\egroup^{\rule{0mm}{2mm}}\right.}

\def\car#1{\mbox{{\rm CAR$\left({#1}^{\vplatz{1.5mm}}\right)$}}}
\def\ccr#1{\mbox{{\rm CCR$\left({#1}^{\vplatz{1.5mm}}\right)$}}}

  \def\ccr #1,#2.{\overline{\Delta(#1,\,#2)}}

  \def\b #1.{{\bf #1}}
  \def\cross#1.{\mathrel{\mathop{\times}\limits_{#1}}}
  \def\C{\Mb{C}}
  \def\N{\Mb{N}}
  \def\R{\Mb{R}}
  \def\Z{\Mb{Z}}

  \def\wh{\widehat}
  \def\wwh #1.{\widehat{#1}}
  \def\wt #1.{\widetilde{#1}}

  \def\cross #1.{\mathrel{\raise 3pt\hbox{$\mathop\times\limits_{#1}$}}}
\def\set #1,#2.{\left\{\,#1\;\bigm|\;#2\,\right\}}
\def\b #1.{{\bf #1}}

\def\ol #1.{\overline{#1}}
\def\rn#1.{\romannumeral{#1}}

\def\rest{\restriction}

\def\s #1.{_{\smash{\lower2pt\hbox{\mathsurround=0pt $\scriptstyle #1$}}\mathsurround=3pt}}
\def\bra #1,#2.{{\left\langle #1,\,#2\right\rangle_{\al A.}}}

\def\XP#1!{\renewcommand{\baselinestretch}{.7}\marginpar{{\footnotesize #1}\hfil}
\renewcommand{\baselinestretch}{1.5}}
\def\XB{\marginpar{
{\footnotesize\bf Change~starts----}\lower 11pt\hbox{\mathsurround=0pt$
\!\!\displaystyle{
\Bigg\downarrow}$\mathsurround=3pt}}}
\def\XE{\marginpar{{\footnotesize\bf Change~ends-----}\raise 10pt\hbox{\mathsurround=0pt$
\!\!\displaystyle{
\Bigg\downarrow}$\mathsurround=3pt}}}

\newcommand{\PP}{\mathbb{P}}

\newcommand{\cB}{{\Sigma}}

\newcommand{\cH}{{\mathcal H}}

\newcommand{\cK}{{\mathcal K}}
\newcommand{\cL}{{\mathcal L}}
\newcommand{\cM}{{\mathcal M}}
\newcommand{\cN}{{\mathcal N}}

\title{Modular Theory by example}
\author{Fernando Lled\'o}
\address{Department of Mathematics, University Carlos~III Madrid,
     Avda.~de la Universidad~30, E-28911 Legan\'es (Madrid),
     Spain and Institute for Pure and Applied Mathematics,
     RWTH-Aachen University,
     Templergraben 55, D-52062 Aachen, Germany (on leave)}
\email{flledo@math.uc3m.es {\em and} lledo@iram.rwth-aachen.de}

\date{\today{}}
\begin{document}
\maketitle

\tableofcontents

\begin{abstract} 
  The present article contains a short introduction to Modular Theory
  for von Neumann algebras with a cyclic and separating vector. It
  includes the formulation of the central result in this area, the
  Tomita-Takesaki theorem, and several of its consequences. We
  illustrate this theory through several elementary examples. We also
  present more elaborate examples and compute modular objects for a
  discrete crossed product and for the algebra of canonical
  anticommutation relations (CAR-algebra) in a Fock representation.
\end{abstract}

\section{Introduction}

Modular Theory has been one of the most
exciting subjects for operator algebras and for its
applications to mathematical physics. We will give here
a short introduction to this theory and state some of its main results.
There are excellent textbooks and review articles which 
cover this subject, e.g.~\cite{bStratila81,bTakesakiII,bSunder87,DaeleIn76,
bFillmore96}, \cite[Section~2.5.2]{bBratteli87} or \cite[Chapter~9]{bKadisonII}.
For an overview and further applications to quantum field
theory see also \cite{Borchers00,pBorchers04,SummersIn04} and references 
cited therein.
A beautiful alternative approach to Modular Theory in terms of bounded operators
is given in \cite{Rieffel77}. This approach is close in spirit to the example
presented in the context of the CAR-algebra in Section~\ref{sec:CAR}.
The origin of the terminology is explained in Example~\ref{OriginModular}.

This article is not intended as a systematic study of Modular Theory
for von Neumann algebras. Rather, the
emphasis lies on the examples. The hope is that the reader will
recognize through the 
examples some of the power, beauty and variety of 
applications of Modular Theory.
We have also included a few exercises to motivate further thoughts
on this topic.
Additional aspects and applications of this theory 
will also appear in \cite{pGuido08}.

In the present article we will present Modular Theory in the special 
case when the von Neumann algebra $\al M.$ and its commutant $\al M.'$
have a common cyclic vector $\Omega$. This approach avoids introducing too much
notation and is enough for almost all examples and applications presented
in this school. (An exception to this is Example~\ref{OriginModular}.)
The reader interested in the more general context described in terms
of Hilbert algebras is referred to \cite[Chapter~VI]{bTakesakiII} or 
\cite{bTakesaki70,bStratila81}.
\section{Modular Theory: definitions, results and first 
examples}\label{MatModularTheory}

In Modular Theory one studies systematically the relation of a von Neumann
algebra $\al M.$ and its commutant $\al M.'$ in the case where both algebras
have a common cyclic vector $\Omega$. We begin introducing some standard 
terminology and stating some elementary results: 

\begin{definition}
Let $\al M.$ be a von Neumann algebra on a Hilbert space $\ot h.$. 
A vector $\Omega\in\ot h.$ is called {\em cyclic for $\al M.$} 
if the set $\{M\Omega\mid M\in\al M.\}$ is dense in $\ot h.$.
We say that $\Omega\in\ot h.$ is {\em separating for $\al M.$}
if for any $M\in\al M.$, $M\Omega=0$ implies $M=0$.
\end{definition}

\begin{proposition}\label{CycSep}
Let $\al M.$ be a von Neumann algebra on a Hilbert space $\ot h.$ and
$\Omega\in\ot h.$. Then
$\Omega$ is cyclic for $\al M.$ iff
$\Omega$ is separating for $\al M.'$.
\end{proposition}
\begin{proof}
Assume that $\Omega$ is cyclic for $\al M.$ and take $M'\in\al M.'$
such that $M'\Omega=0$. Then $M'M\Omega=MM'\Omega=0$ for all $M\in\al M.$.
Since $\{M\Omega\mid M\in\al M.\}$ is dense in $\ot h.$ it follows that
$M'=0$.

Assume that $\Omega$ is separating for $\al M.$ and denote by $P'$
the orthogonal projection onto the closed subspace generated by
$\{M\Omega\mid M\in\al M.\}$. Then $P'\in\al M.'$ and 
$(1-P')\Omega=\Omega-\Omega=0$. Since $\Omega$ is separating 
for $\al M.'$ we 
have that $P'=\1$, hence $\Omega$ is also cyclic for $\al M.$.
\end{proof}

If $\Omega\in\ot h.$ is cyclic for the von Neumann algebras $\al M.$
and its commutant $\al M.'$ (hence also separating for both algebras
by the preceding proposition) one can naturally introduce the
following two antilinear operators $S_0$ and $F_0$ on $\ot h.$:
\begin{eqnarray*}
  S_0( M\Omega)&:=&M^*\,\Omega\;,\quad M\in\al M. \\
  F_0( M'\Omega)&:=&(M')^*\,\Omega\;,\quad M'\in\al M.'\,.
\end{eqnarray*}
Both operators are well defined on the dense domains
$\mr dom.S_0=\al M.\Omega$ and $\mr dom.F_0=\al M.'\Omega$,
respectively, and have dense images.
It can be shown that the operators $S_0$ and $F_0$ are
closable and that $S=F_0^*$ as well as $F=S_0^*$,
where $S$ and $F$ denote the closures of $S_0$ and $F_0$, respectively.
The closed, antilinear operator $S$ is called the {\em Tomita operator}
for the pair $(\al M.,\Omega)$, where $\Omega$ is cyclic and separating
for $\al M.$. The operators $S$ and $F$ are involutions in the sense that
if $\xi\in\mr dom.S$, then $S\xi\in\mr dom.S$ and $S^2\xi=\xi$
(similarly for $F$).

Let $\Delta$ be the unique positive, selfadjoint operator and
$J$ the unique antiunitary operator
occurring in the polar decomposition of $S$, i.e.
\[
 S= J\,\Delta^\frac12\,.
\] 
We call $\Delta$ the {\em modular operator}
and $J$ the {\em modular conjugation}
associated with the pair $(\al M.,\Omega)$.

We mention next standard relations between the previously defined modular objects 
$S,F,\Delta$ and $J$. For a complete proof see Proposition~2.5.11 in \cite{bBratteli87}.
\begin{proposition}\label{properties}
The following relations hold
\begin{eqnarray*}
  \Delta &=& FS\;,\quad \Delta^{-1}\;\;=\;\; SF\;,
             \quad F\;\;=\;\; J\Delta^{-\frac12}\\
  J &=& J^* \;,\quad J^2\;\;=\;\; \1\;,
        \quad \Delta^{-\frac12} \;\;=\;\; J\Delta^{\frac12} J\,.
\end{eqnarray*}
\end{proposition}

We conclude stating the main result of modular
theory, the so-called Tomita-Takesaki theorem. We will give a 
proof only for the case where all modular objects are bounded.
This situation covers Examples~\ref{ex:LeftReg} 
and \ref{ExFinBg}, as well as the results in Section~\ref{sec:ex-crossed}
below. To state the theorem we need to introduce the following notation:
given the modular operator $\Delta$,
we construct the strongly continuous unitary group
\[
 \Delta^{it}=\mr exp.\Big(\,it \left(\mr ln.\Delta\right) \Big)
             \;,\quad t\in\R\,,
\]
via the functional calculus.
It is called the {\em modular group} and 
\[
 \sigma_t(M):=\Delta^{it}\,M\,\Delta^{-it}\;,\quad M\in\al M.\,,\,t\in\R
\]
gives a one parameter automorphism group on $\al M.$, the so-called
{\em modular automorphism group}.

\begin{theorem} (Tomita-Takesaki)\label{TTT}
Let $\al M.$ be a von Neumann algebra with cyclic 
and separating vector $\Omega$ in the Hilbert space
$\ot h.$. The operators $\Delta$ and $J$ are the 
corresponding modular operator and 
modular conjugation, respectively, and denote by $\sigma_t$, $t\in\R$,
the modular automorphism group. Then we have
\[
 J\,\al M.\, J=\al M.'\quad\mr and.\quad 
 \sigma_t(\al M.)=\al M.\;,\;\; t\in\R\,.
\]
\end{theorem}
\begin{proof}
We will proof\footnote{This proof is close to the one given in
\cite[p.~48-49]{bSunder87}. See also \cite{DaeleIn82}.} 
this result assuming that the Tomita $S$ operator
is bounded. This implies that $S^*$, $\Delta=S^*S$ and $\Delta^{-1}=SS^*$
are also bounded.

i) We show first $S\al M.S=\al M.'$. For any $M,M_0\in\al M.$ and using
$S\Omega=\Omega$ we have
\begin{equation}\label{change}
 SM_0S\,M\Omega=SM_0M^*\Omega=MM_0^*\Omega=M\,S(M_0\Omega)=M\,S M_0 S\Omega
\end{equation}
or, equivalently, $(SM_0S\,M\Omega-M\,S M_0 S)\Omega=0$. Putting $M=M_1M_2$
and using again Eq.~(\ref{change}) we get
\begin{eqnarray*}
(SM_0S\,M_1M_2-M_1M_2\,S M_0 S)\;\Omega&=&0\quad\mathrm{for~all}\; M_2\in\al M.\\
(SM_0S\,M_1-M_1\,S M_0 S)\;M_2\,\Omega&=&0\quad\mathrm{for~all}\; M_2\in\al M.\,.
\end{eqnarray*}
Since $\{M\Omega\mid M\in\al M.\}$ is dense in $\ot h.$ it follows that
$S\al M.S\subseteq\al M.'$. Similarly it can be shown that 
$S^*\al M.'S^*\subseteq\al M.''=\al M.$. Taking in the last 
inclusion adjoints and 
multiplying both sides with $S$ we get the reverse inclusion 
$\al M.'\subseteq S\al M.S$, hence
\[
S\al M.S=\al M.'\quad\mathrm{and} \quad S^*\al M.'S^*=\al M.\,.
\]

ii) Next we prove the following statement: $\Delta^z\al
M.\Delta^{-z}=\al M.$, $z\in\C$. It implies the required equation
taking $z=it$. By the preceding item and using 
$\Delta=S^*S$, $\Delta^{-1}=SS^*$ we get
\begin{equation}\label{cero}
 \Delta\al M.\Delta^{-1}=S^*S\,\al M.\,SS^*=\al M.\,,\quad
 \mathrm{hence}\quad\Delta^n\al M.\Delta^{-n}=\al M.\;,\;n\in\N\,.
\end{equation}
Since $\Delta$ is bounded, selfadjoint and invertible we may apply 
the spectral theorem and functional calculus to obtain the following
representations
\[
 \Delta=\int_{\varepsilon}^{\|\Delta\|} \lambda\,dE(\lambda)\quad
 \mathrm{and}\quad
 \Delta^z=\int_{\varepsilon}^{\|\Delta\|} \lambda^z\,dE(\lambda)
 \;,\quad z\in\C\,,
\]
for some positive $\varepsilon$. It can be shown
that for any $M\in\al M.$, $M'\in\al M.'$ and
$\xi,\eta\in\ot h.$ the function defined by
\[
 f(z)=\|\Delta\|^{-2z}
      \left\langle\left[\Delta^z M\Delta^{-z},M'\right]\xi,\eta\right\rangle
\]
is entire. Moreover, it is also bounded since
\[
|f(z)|\leq \|\Delta\|^{-2\mathrm{Re}(z)}\,2\,\|\Delta^z\|\|\Delta^{-z}\|
           \|M\|\|M'\|\|\xi\|\|\eta\|
      \leq 2\|M\|\|M'\|\|\xi\|\|\eta\|\;,
\]
where the last inequality holds if Re$(z)\geq 0$ and we have used
$\|\Delta^z\|=\|\Delta^{-z}\|$. 
Altogether, we have
constructed an entire function $f$ which is bounded if Re$(z)\geq 0$ and 
by Eq.~(\ref{cero}) satisfies $f(n)=0$, $n\in\N$. Therefore $f(z)=0$ for
all $z\in\C$ and
\begin{equation}\label{inclusion-z}
\Delta^z\al M.\Delta^{-z}\subseteq\al M.''=\al M.\;,\quad z\in\C\;.
\end{equation}
Multiplying the last inclusion from the left by $\Delta^{-z}$ and
from the right by $\Delta^z$ and using (\ref{inclusion-z}) changing $z$ 
by $-z$ we obtain $\Delta^z\al M.\Delta^{-z}=\al M.$.

iii) Finally we show the relation involving the modular conjugation. Using 
the preceding step (ii) we get
\[
 J\al M.J =J\Delta^\frac12\,\al M.\,\Delta^{-\frac12} J=S\al M.S^*=\al M.'
\]
and the proof is concluded.
\end{proof}

\subsection{Comments and elementary consequences of the Tomita-Takesaki 
            theorem}

One can recognize in the first part of the preceding theorem the 
interplay between algebraic and analytic structures 
(cf.~\cite{pLledo08a}). In fact, the commutant of a von Neumann
algebra is obtained by conjugation with an {\em analytic} object like $J$,
which is obtained in terms of the polar decomposition 
of an antilinear closed operator.

\begin{remark}\label{rem-mod}
\begin{itemize}
\item[(i)] There are various approaches to a 
 complete  proof of Theorem~\ref{TTT}:
 one of them stresses more the analytic aspects and techniques from
 the theory of unbounded operators; another one 
 emphasizes more the algebraic structure (cf.~\cite[Chapter~9]{bKadisonII}).
 In \cite{Rieffel77} Rieffel and van Daele present a different proof based on projection 
 techniques and bounded operators. This approach is justified by the fact
 that the main ingredients of the theorem, namely the modular conjugation 
 $J$ and the modular group $\Delta^{it}$ can be characterized in terms
 of real subspaces which have suitable relative positions
 within the underlying complex Hilbert space.
\item[(ii)] An immediate application of the preceding theorem is that the 
modular conjugation $J\colon\ot h.\to\ot h.$ is a *-anti-isomorphism between 
$\al M.$ and its commutant $\al M.'$. 
\item[(iii)] Assume that the cyclic and separating vector $\Omega$ 
for the von Neumann algebra has norm $1$. 
Then the Tomita operator $S$ measures to what extent the corresponding
vector state $\omega$ on $\al M.$ defined by
\begin{equation}\label{state}
\omega(M):=\langle\Omega,M\Omega\rangle \;, \quad M\in\al M.
\end{equation}
is tracial.
In fact, note that $S$ is an isometry iff $\omega$ is a trace, since
\[
\|M\Omega\|^2=\omega(M^*M)=\omega(MM^*)=\|M^*\Omega\|^2=\|S(M\Omega)\|^2\,.
\]
The vector state $\omega$ associated to $\Omega$ is a faithful normal state.
Conversely, to any faithful normal state of $\al M.$ one can associate, via
the GNS construction, a cyclic and separating vector in the GNS Hilbert space.
Modular Theory may be extended to the situation of von Neumann algebras with
faithful, normal and semifinite weights (see e.g.~\cite{bSunder87}).
\end{itemize}
\end{remark}

The following proposition, which can be shown directly (recall the exercise proposed
in \cite[Subsection~2.2]{pLledo08a}), 
is an easy to prove if we use Theorem~\ref{TTT}.
\begin{proposition}\label{Pro:MaxAb}
Let $\al A.\subset\al L.(\ot h.)$ be an Abelian von Neumann algebra with a 
cyclic vector $\Omega\in\ot h.$. Then $\al A.$ is maximal Abelian, 
i.e.~$\al A.=\al A.'$.
\end{proposition}
\begin{proof}
Since the algebra $\al A.$ is Abelian we have $\al A.\subseteq\al A.'$ and 
any cyclic vector $\Omega$ for $\al A.$ 
will also be cyclic for $\al A.'$. Therefore we can apply Theorem~\ref{TTT}
and the following chain of inclusions
\[
 \al A.\subseteq\al A.'= J\,\al A.\,J\subseteq J\,\al A.'\,J 
                         \subseteq J\,(J\,\al A.\,J)\,J 
                       = \al A.\,.
\]
conclude the proof.
\end{proof}

\begin{remark}
One of the origins of Modular Theory can be traced back to the original
work by Murray and von Neumann. 
A vector $u\in\ot h.$ is called a trace
vector for a von Neumann algebra $\al M.\subset\al L.(\ot h.)$ if
\[
 \langle u,MN u\rangle=\langle u,NM u\rangle\;,\quad M,N\in\al M.\;.
\]
If $\al M.$ has a cyclic trace vector $u$, then for any $M\in\al M.$ 
there is a unique $M'\in\al M.'$ satisfying
\[
 Mu=M'u\;.
\]
In this case we say that $M$ and $M'$ are reflections of one another
about $u$. The mapping $M\to M'$ is a *-anti-isomorphism between $\al
M.$ and $\al M.'$. Nevertheless, this result remains of limited
applicability since the assumption that the von Neumann algebra has a
generating trace vector is so strong that it already implies that $\al
M.$ and $\al M.'$ are finite (see \cite[Theorem~7.2.15]{bKadisonII} or
\cite{KadisonIn91} for more details).
\end{remark}

\subsection{Examples}

The construction of the modular objects given in the
beginning of this section is rather involved. The modular
operator $\Delta$ and modular conjugation $J$ appear as 
the components of the polar decomposition associated to the 
closure of the operator $S_0$ defined above. To gain some
intuition on the modular objects it is useful to compute 
$J$ and $\Delta$ in concrete cases. We begin with some 
natural examples related to the representation theory of groups.

\begin{example}\label{ex:LeftReg}
We will see in this example that the corresponding modular 
objects are bounded.
Let $\al G.$ be a discrete group and consider the Hilbert space
$\al H.:=\ell_2(\al G.)$ with orthonormal basis given by delta functions
on the group $\{\delta_g(\cdot)\mid g\in\al G.\}$. Define the 
left- resp.~right regular (unitary) representations on $\al H.$ as
\begin{equation}\label{def:LR}
 L(g_0)\delta_g:=\delta_{g_0g}\quad\mr resp..\quad
 R(g_0)\delta_g:=\delta_{gg_0^{-1}}\;,\quad g\in\al G.\,.
\end{equation}
We introduce finally the von Neumann algebra generated by the left
regular representation
\[
 \al M.:=\{L(g)\mid g\in\al G. \}''\subset\al L.(\al H.)\,.
\]
Using Eq.~(\ref{def:LR}) it
is immediate to verify that $\{R(g)\mid g\in\al G. \}\subset\al M.'$
(hence $\{R(g)\mid g\in\al G. \}''\subseteq\al M.'$)
and that $\Omega:=\delta_e$ is a cyclic vector
for $\al M.$ and $\al M.'$. To determine the Tomita operator it 
is enough to specify $S$ on the canonical basis of $\ell_2$ and extend this 
action anti-linearly on the whole Hilbert space:
\[
 S(\delta_g)=S\Big(L(g)\delta_e\Big)=L(g)^*\delta_e
            =L(g^{-1})\delta_e = \delta_{g^{-1}}\,.
\]
This implies
\[
  S^*=S\,,\quad \Delta=\1\quad\mathrm{and}\quad J=S\,.
\]
Moreover, it is also straightforward to check the relation between the left- and 
right regular representations in terms of $J$:
\[
  J L(g) J = R(g)\;,\quad g\in\al G.\;.
\]
Finally, we can improve the inclusion
$\{R(g)\mid g\in\al G. \}''\subseteq\al M.'$ mentioned above
applying the first part of Theorem~\ref{TTT}. In fact,
the commutant of $\al M.$ is generated, precisely, by the right regular 
representation:
\[
 \al M.'= J\al M. J=J\{L(g)\mid g\in\al G. \}''J
                   =\{JL(g)J\mid g\in\al G. \}''
                   =\{R(g)\mid g\in\al G. \}''\,.
\]
\end{example}

\begin{example}\label{OriginModular}
Apparently one of the original motivations of Tomita for developing 
Modular Theory was the harmonic analysis of nonunimodular locally
compact groups. For the following example the notion of Hilbert algebras is needed
(see \cite[Sections~2.3 and 2.4]{bSunder87} for details).
Let $\al G.$ be a locally compact group with left invariant 
Haar measure $dg$ and modular function
\[     
\widetilde{\Delta}(\cdot): \al G.\to \R_+\,.
\]
(Recall that the modular function is a continuous group homomorphism
that relates the left and right Haar integrals, 
cf.~\cite[\S 15]{bHewittI}). As in the preceding example,
the modular objects associated to the left regular representation 
on the Hilbert space $L^2(\al G.,dg)$ are given in this case by
\begin{itemize}
\item[] $(S\varphi)(g)=\widetilde{\Delta}(g)^{-1}\;\overline{\varphi(g^{-1})}$,
        $\qquad \varphi\in L^2(\al G.,dg)$
\item[] $(J\varphi)(g)=\widetilde{\Delta}(g)^{-\frac12}\;\overline{\varphi(g^{-1})}$,
        $\qquad \varphi\in L^2(\al G.,dg)$
\item[] $(\Delta\varphi)(g)=\widetilde{\Delta}(g)\;\varphi(g)$, 
        $\qquad\qquad \varphi\in L^2(\al G.,dg)$.
\end{itemize}
This example shows the origin of the name Modular Theory,
since the modular operator is just multiplication by the 
modular function of the group. Moreover, the preceding
expressions of the modular objects
are in accordance with the preceding example. Recall
that if $\al G.$ is discrete, then it is also unimodular, 
i.e.~$\Delta(g)=1$, $g\in\al G.$.
\end{example}

\begin{example}\label{ExFinBg}
Let $\al H.$, $\al H.'$ be finite dimensional Hilbert spaces with
$\mathrm{dim}\al H.=\mathrm{dim}\al H.'=n$ and orthonormal 
basis $\{e_k\}_{k=1}^n$ and $\{e_k'\}_{k=1}^n$, respectively.
Consider on the tensor product $\al H.\otimes\al H.'$ the von Neumann
algebra $\al M.=\al L.(\al H.)\otimes\C\1_{\al H.'}$. It is easy to verify
that the vector
\[
 \Omega:=\sum_k \lambda_k (e_k\otimes e_k')\in\al H.\otimes\al H.'\,,\;
 \mbox{with}\; \lambda_k>0\,,\;k=1,\ldots,n\,,\; \mbox{and}\;
 \sum_{k=1}^n \lambda_k^2=1\,,
\] 
is a cyclic and separating vector for $\al M.$ with norm 1.

A direct computation shows that 
the modular objects for the pair $(\al M.,\Omega)$ are given by
\begin{eqnarray*}
S(e_k\otimes e_s')  &=&\frac{\lambda_k}{\lambda_s}\,(e_s\otimes e_k')\,,\quad
                          k,s\in\{1,\dots,n\}.\\
S^*(e_k\otimes e_s')&=&\frac{\lambda_s}{\lambda_k}\,(e_s\otimes e_k')\,,\quad
                          k,s\in\{1,\dots,n\}.\\
\Delta(e_k\otimes e_s')
                    &=&\left(\frac{\lambda_k}{\lambda_s}\right)^2\,(e_k\otimes e_s').\\
J(e_k\otimes e_s')  &=&(e_s\otimes e_k').
\end{eqnarray*}
In this example the modular conjugation $J$ acts as a flip of indices
for the given basis of the tensor product Hilbert space. This action
extends anti-linearly to the whole tensor product $\al H.\otimes\al
H.'$. Using Theorem~\ref{TTT} and the explicit expression for $J$ we
can again improve the inclusion $\C\1\otimes\al L.(\al H.)\subseteq\al M.'$.
In fact, note that for any $A\in\al L.(\al H.)$ we have
\[
 J\,(A\otimes \1)\,J=(\1\otimes \overline{A})\;,\;\mathrm{where}\;\; 
                     (\overline A)_{rs}=(\overline{A_{rs}})\;.
\]
Now applying Theorem~\ref{TTT} we get
\[
 \al M.'= J\,\al M.\, J = J\,(\al L.(\al H.)\otimes\C\1)\, J=\C\1\otimes\al L.(\al H.)\;.
\]
\end{example}

\begin{exercise}
Generalize the preceding example to the case where $\al H.$ is an infinite dimensional
separable Hilbert space.
\end{exercise}

\section{Modular objects for a crossed product}\label{sec:ex-crossed}

The crossed product construction is a procedure to obtain a new von Neumann algebra 
out of a given von Neumann algebra which carries a certain group action. This
principle, in particular the group measure space construction given below,
goes back to the pioneering work of Murray and von Neumann on rings of 
operators (cf.~\cite{vNeumannI,vNeumannIII}). 
Standard references which present the crossed product construction
with some variations are
\cite[Chapter~4]{bSunder87}, \cite[Section~V.7]{bTakesakiI} or
\cite[Section~8.6 and Chapter~13]{bKadisonII}. 

Let $(\Omega,\Sigma,\PP)$ be a separable measure space with
probability measure $\PP$ defined on Borel $\sigma$-algebra $\cB$.
Consider the Hilbert space $\cH:=L^2(\Omega,\PP)$ and identify with
$\cM:=L^\infty(\Omega,\PP)$ the abelian von Neumann algebra that acts
on $\cH$ by multiplication.  Suppose that there is an infinite
countable discrete group $\Gamma $ acting on $(\Omega,\cB,\PP)$ by
measure preserving automorphisms\footnote{An automorphism $T$ of the measure
space $(\Omega,\Sigma,\PP)$ is a bijection $T\colon\Omega\to\Omega$
such that (i) for $\al S.\in\Sigma$ we have $T(\al S.),T^{-1}(\al S.)\in\Sigma$
and (ii) if $\al S.\in\Sigma$, then $\mu(\al S.)=0\Leftrightarrow\mu(T^{-1})=0$.
The action $T$ is measure preserving if $\mu\circ T_\gamma=\mu$ for
all $\gamma\in\Gamma$.}, 
i.e., $T\colon \Gamma
\to\mathrm{Aut}(\Omega,\cB,\PP)$. This action induces a canonical
action $\alpha$ of $\Gamma $ on the von Neumann algebra $\cM$:
\begin{equation*}
\alpha\colon \Gamma \to\mathrm{Aut}\cM\quad\mathrm{with}\quad (\alpha_\gamma
f)(\omega):=(f\circ T_\gamma^{-1})(\omega)\;,\quad f\in\al M., \omega\in\Omega\,.
\end{equation*}
\begin{definition}
Let $\alpha\colon \Gamma \to\mathrm{Aut}\cM$ be an action of the discrete group
$\Gamma$ on $\al M.$ as above. 
  \begin{itemize}
  \item[(i)] The action $\alpha$ is called free if any
    $\alpha_\gamma$, $\gamma\not=e$, satisfies the following
    implication: the equation $fg=\alpha_\gamma(g)f$ for all $g\in\al
    M.$ implies $f=0$.
  \item[(ii)] The action $\alpha$ is called {\em ergodic} if the
    corresponding fixed point algebra is trivial, i.e.
\[
 \al M.^\alpha:=\{f\in\al M.\mid \alpha_\gamma(f)=f \;,\quad\gamma\in\Gamma \}=\C 1\;.
\]
  \end{itemize}
  
\end{definition}

\begin{remark}
  The preceding properties of a group action on the von Neumann
  algebra also translate into properties of the corresponding group
  action $T\colon\Gamma \to\mathrm{Aut}(\Omega,\cB,\PP)$ on the
  probability space.
 \begin{itemize}
  \item[(i)] The action $T$ is called {\em free} if the sets $\{\omega\in\Omega\mid
  T_\gamma\omega=\omega\}$ are of measure zero for all $\gamma\not=e$.
  \item[(ii)] The action $T$ is called ergodic if for any $\al S.\in\Sigma$ such that
\[
   \mu\left(\left(T_\gamma(\al S.)\setminus\al S.\right)\cup
   \left(\al S.\setminus T_\gamma(\al S.)\right)\right)=0\;,\quad\mathrm{for~all} \quad
   \gamma\in\Gamma\;,
\] 
one has either $\mu(\al S.)=0$ or $\mu(X\setminus \al S.)=0$.
 \end{itemize}
\end{remark}

The crossed product is a new von Neumann algebra $\al N.$ which can be
constructed from the dynamical system $(\al M.,\alpha,\Gamma)$. It
contains a copy of $\al M.$ and a copy of the discrete group $\Gamma$ by
unitary elements {\em in} $\al N.$ and the commutation relations
between both ingredients are given by the group action.  We will
describe this construction within the group measure space context
specified above.  We begin introducing a new Hilbert space on which
the crossed product will act:
\begin{equation}\label{eq:decomp}
 \cK:=\ell_2(\Gamma )\otimes\cH\cong\bigoplus_{\gamma\in\Gamma }\cH_\gamma\,
                           \cong \int_\Omega \ell^2(\Gamma ) d\PP,
\end{equation}
where $\cH_\gamma\equiv\cH=L^2(\Omega,\PP)$ for all $\gamma\in\Gamma$. 
Moreover, we consider the following representations of 
$\al M.=L^\infty(\Omega,\PP)$ and $\Gamma$
on $\al K.$: for $\xi=(\xi_\gamma)_{\gamma\in\Gamma}\in\al K.$ with
$\xi_\gamma\in\al H.$ we define
\begin{eqnarray}\label{abelianM}
(\pi(f)\xi)_\gamma&:=&\alpha^{-1}_\gamma(f)\xi_\gamma
= f(T_\gamma \cdot)\xi_\gamma\;,\quad f\in\cM, \\ \label{lrr}
(U(\gamma_0)\xi)_\gamma&:=&\xi_{\gamma_0^{-1}\gamma} \,.
\end{eqnarray}

The {\em discrete crossed product} of $\cM$ by $\Gamma $ is the von
Neumann algebra acting on $\cK$ and generated by these operators, i.e.,
\begin{equation*}
  \cN=\cM\otimes_\alpha\Gamma 
    := \Big(
        \left\{\pi(f)\mid f\in\al M.\right\}\cup
        \left\{ U(\gamma)\mid\gamma\in\Gamma \right\}
       \Big)''
     \subset\al L.(\cK)\,,
\end{equation*} where the prime denotes the commutant in
$\cL(\cK)$. A characteristic relation for the crossed product is
\begin{equation}\label{main}
\pi\left(\alpha_\gamma(f)\right)=U(\gamma)\pi(f)U(\gamma)^{-1}\,.
\end{equation}
In other words, $\pi$ is a covariant representation of the $W^*$-dynamical system
$(\cM, \Gamma ,\alpha)$. 

\begin{remark}\label{circulant}
   It is useful to characterize explicitly {\em all} elements in the 
   crossed product and not just a generating family. For this consider
   the identification $\cK\cong\bigoplus_{\gamma\in\Gamma }\cH_\gamma$ with
   $\cH_\gamma\equiv\cH$. Then every $T\in\cL(\cK)$ can be written as an infinite
   matrix $(T_{\gamma'\gamma})_{\gamma',\gamma\in\Gamma }$ with entries
   $T_{\gamma'\gamma}\in\cL(\cH)$. Any element $N$ in the crossed
   product $\cN\subset\cL(\cK)$ has the form
\begin{equation}\label{eq:matrix-n}
   N_{\gamma'\gamma}=\alpha_\gamma^{-1}\left( A(\gamma'\gamma^{-1})\right)
                   \;,\quad \gamma'\,,\,\gamma\in\Gamma \,,
 \end{equation}
 for some function $A\colon\Gamma \to\cM\subset\cL(\cH)$. 
 Since any $N\in\cN$ is a bounded operator in $\al K.$ it follows that
\begin{equation}\label{series-condition}
 \sum_{\gamma\in\Gamma} \|A(\gamma) \|^2_\al H. < \infty\;.
\end{equation}
  
  For example, the matrix expression of the product of generators 
  $N:=\pi(g)\cdot U(\gamma_0)$, $g\in\cM$, $\gamma_0\in\Gamma$, is given by
\begin{equation}\label{matrix-generators}
 N_{\gamma'\gamma}=\alpha_{\gamma}^{-1}(g) \, \delta_{\gamma',\gamma_0\gamma}
                  =\alpha_{\gamma}^{-1}\left(A(\gamma'\gamma^{-1})\right)\;,
\end{equation}
where
$ A(\widetilde{\gamma}):=\left\{\begin{array}{l}
                               g\quad\mathrm{if}\;\;\widetilde{\gamma}=\gamma_0\\
                               0\quad\mathrm{otherwise}
                               \end{array}\right.$.
\end{remark}

With the group measure space construction one can produce 
examples of von Neumann algebras of any type. In this section 
we concentrate on the case of finite factors. 
(For a general statement see Theorem~\ref{ex-factors}.)
\begin{theorem}\label{teo:finite}
If the action of the discrete group $\Gamma$ on the probability space 
$(\Omega,\Sigma,\PP)$ is measure preserving, free and ergodic, then the crossed product 
$\al N.$ constructed before is a finite factor, i.e. a factor of type I$_n$
or of type II$_1$.
\end{theorem}

Next we specify a cyclic and separating vector for the crossed product $\al N.$
defined previously.
\begin{lemma}\label{lem:cyc-sep}
The vector $\Omega\in\al K.$ defined by $(\Omega)_\gamma=\delta_{e\gamma}1$, where
$1\in\al H.$ is the identity function, is cyclic and separating for
the finite factor $\al N.$. 
\end{lemma}
\begin{proof}
1. We show first that $\Omega$ is separating for $\al N.$. Let $N\Omega=0$ for 
some $N\in\al N.$. Then by the matrix form of the elements of $\al N.$ mentioned
in Remark~\ref{circulant} we have
\[
 0=(N\Omega)_\gamma=A(\gamma)1\;,\quad \mathrm{for~all}\;\gamma\in\Gamma\,.
\]
This implies $A=0$, hence  $N=0$.

2. To show that $\Omega$ is cyclic note that it is enough to verity that the set
\[
 \al D.:=\{\ot g.=(g_\gamma)_\gamma\mid g_\gamma\in\al M.\;\;
           \mathrm{and~finitely~many}\; g_\gamma\not=0\}
\]
is dense in $\al K.$. In fact, using Eq.~(\ref{matrix-generators}) we have
\[
 \al D.=\mathrm{span}\{\pi(g)\cdot U(\gamma)\;\Omega
                       \mid g\in\al M.\;,\; \gamma\in\Gamma\}
        \subset \{N\Omega\mid N\in\al N.\}\,.
\]
Now, for any $\varphi=(\varphi_\gamma)_\gamma\in\al K.$ and any $\varepsilon>0$
choose a subset $\Gamma_0\subset\Gamma$ with finite cardinality, 
i.e.~$|\Gamma_0|<\infty$, such that
\[
  \sum_{\gamma\in (\Gamma_0)^c}\|\varphi_\gamma\|^2_\al H.<\frac{\varepsilon^2}{2}\,.
\]
Since $\al M.=L^\infty(\Omega,\PP)$ is dense in $\al H.=L^2(\Omega,\PP)$ choose
an element $\ot g.=(g_\gamma)_\gamma\in\al D.$ such that $g_\gamma=0$ for 
$\gamma\in(\Gamma_0)^c$ and
\[
\|\varphi_\gamma-g_\gamma\|_\al H.\leq\frac{\varepsilon^2}{2|\Gamma_0|}\quad
 \mathrm{for}\; \gamma\in\Gamma_0\,.
\]
Then we have 
\[
 \|\varphi-\ot g.\|^2_\al K.
              =\sum_{\gamma\in \Gamma_0}\|\varphi_\gamma-g_\gamma\|^2_\al H.+
               \sum_{\gamma\in (\Gamma_0)^c}\|\varphi_\gamma\|^2_\al H.
              \leq |\Gamma_0|\cdot\frac{\varepsilon^2}{2|\Gamma_0|}
                   +\frac{\varepsilon^2}{2}=\varepsilon^2\,.
\]
This shows that $\al D.$ is dense in $\al K.$ and the proof is concluded.
\end{proof}

\begin{remark}
  In the proof of the cyclicity of $\Omega$ in the preceding lemma it
  is crucial that $L^\infty(\Omega,\PP)$ is dense in
  $L^2(\Omega,\PP)$ or, equivalently, $1\in L^2(\Omega,\PP)$ is cyclic
  with respect to $L^\infty(\Omega,\PP)$. This is true because we are
  working with a finite measure space. 
\end{remark}

We finish this section specifying the Modular objects for the pair
$(\al N.,\Omega)$.

\begin{theorem}
The Tomita operator associated to the pair $(\al N.,\Omega)$ 
given in Theorem~\ref{teo:finite} and Lemma~\ref{lem:cyc-sep} is an 
isometry, hence $\Delta=\1$ and $S=J$. An explicit expression for the 
Tomita operator is given on the dense set 
\begin{equation}\label{eq:dense}
 \{N\Omega\mid N\in\al N.\}=\{\ot g.=\left(A(\gamma)\right)_\gamma\mid
               \mathrm{for~some~mapping}\;\,A\colon\Gamma \to\cM\}\subset\al K.
\end{equation}
by
\[
  S_{\gamma'\gamma}=\delta_{\gamma'\gamma^{-1}}\;\alpha_\gamma(C(\cdot))\;,
\]
where $C$ means complex conjugation in $\al M.=L^\infty(\Omega,\PP)$. 
The vector state associated to $\Omega$ and defined by
\[
\omega(M):=\langle\Omega,M\Omega\rangle \;, \quad M\in\al M.\;.
\]
is a trace on the finite factor $\al N.$.
\end{theorem}
\begin{proof}
  Note first that Eq.~(\ref{eq:dense}) is an immediate consequence of
  the characterization of the matrix elements of the crossed product
  given in Eq.~(\ref{eq:matrix-n}),
  $N_{\gamma'\gamma}=\alpha_\gamma^{-1}\left(
    A(\gamma'\gamma^{-1})\right)$, for some function $A\colon\Gamma
  \to\cM$. Next we verify that the matrix elements $S_{\gamma'\gamma}$
  given above correspond to the Tomita operator on the subset
  $\{N\Omega\mid N\in\al N.\}$ (which is dense in $\al K.$ by
  Lemma~\ref{lem:cyc-sep}):
\begin{eqnarray*}
 \left(S(N\Omega)\right)_{\gamma}
       &=&\sum_{\gamma'} S_{\gamma\gamma'}\,A(\gamma')
          \;=\;\sum_{\gamma'} \delta_{\gamma(\gamma')^{-1}}\,
                     \alpha_{\gamma'} \left(\overline{A(\gamma')}\right)
          \;=\;\alpha_{\gamma}^{-1}\left(\overline{A(\gamma^{-1})}\right)\\ 
      &=& N^*_{\gamma'\gamma}\,\delta_{\gamma' e}1
          \;=\; \sum_\gamma (N^*)_{\gamma\gamma'}\Omega_{\gamma'}
          \;=\; \left(N^*\Omega\right)_{\gamma}\;.
\end{eqnarray*}
Moreover, for any $\xi=N\Omega=\left(A(\gamma)\right)_\gamma$
we have
\[
 \|S\xi\|^2_{\al K.}=\sum_{\gamma}\|(S\xi)_\gamma\|^2_\al H.
                    =\sum_{\gamma}\|\alpha_{\gamma}^{-1}
                     \left(\overline{A(\gamma^{-1})}\right)\|^2_\al H.
                    =\sum_{\gamma}\|\left({A(\gamma)}\right)\|^2_\al H. 
                    =\|\xi\|^2_{\al K.}\,.
\]
This shows that the Tomita operator $S$ is an isometry on the 
dense subspace $\{N\Omega\mid N\in\al N.\}$, hence it
extends uniquely to an isometry on the whole Hilbert space $\al K.$.
Therefore $\Delta=S^*S=\1$ and $S=J$.
By Remark~\ref{rem-mod}~(iii) it follows that the vector state associated
to the cyclic and separating vector $\Omega$ is a trace.
\end{proof}
\begin{remark}
For an expression of the modular operator on general crossed products which are 
not necessarily finite factors see \cite[{\S}4.2]{bSunder87}.
\end{remark}

\section{Modular objects for the CAR-algebra}\label{sec:CAR}

In this section we construct the modular objects for 
the algebra of the canonical anticommutation relations (CAR-algebra).
This is a more elaborate example that uses standard results
on the CAR-algebra and its irreducible representations (Fock representations).
We give a short review of these results in Subsection~\ref{CARalgebras}.
Let $(\ot h.,\Gamma)$ be a reference space
as in Theorem~\ref{Teo.2.2.9} 
and denote by $\ot q.$ a closed
$\Gamma$-invariant subspace of $\ot h.$. The orthogonal projection 
associated with $\ot q.$ is denoted by $Q$.
We can naturally associate
with the subspace $\ot q.$ a von Neumann 
algebra that acts on the 
antisymmetric Fock space $\ot F.$ characterized by the basis projection $P$
(cf.~(\ref{FockSpace})):
\begin{equation}\label{VNAm}
 \al M.(\ot q.):= \Big( \{ a(q)\mid q\in \ot q.\}\Big)''
                  \subset \al L.(\ot F.)\,.
\end{equation}
In the present subsection we will analyze the modular objects corresponding
to the pair $(\al M.(\ot q.),\Omega)$, where $\Omega$ is the 
so-called Fock vacuum in $\ot F.$.
More details and applications of the Modular Theory in 
the context of the CAR-algebra can be found in 
\cite{Lledo02a} and \cite[Part~I]{hLledo05}.

\subsection{Two subspaces in generic position}

We will address next the question when the Fock vacuum $\Omega$ is 
cyclic and separating for the von Neumann algebra $\al M.(\ot q.)$.
The answer to this question has to do with the relative position that
the subspace $\ot q.$ has with the one particle space $\ot p.:=P\ot h.$.

\begin{proposition}\label{Teo.3.2}
Let $\ot q.$ be a closed $\Gamma$-invariant subspace of $\ot h.$ as
before and denote the one-particle Hilbert space associated to the basis
projection $P$ by $\ot p.$. Then we have:
\begin{itemize}
\item[(i)] The vacuum vector $\Omega$ is cyclic for $\al M.(\ot q.)\;$
  iff $\;\ot p.\cap\ot q.^\perp=\{0\}$. 
 \item[(ii)] The vacuum vector $\Omega$ is separating for $\al M.(\ot q.)$
  iff $\;\ot p.\cap\ot q.=\{0\}$. 
\end{itemize}
\end{proposition}
\begin{proof}
We will only proof part~(i). Similar arguments can be used to show (ii).

Assume that $\Omega$ is cyclic for $\al M.(\ot q.)$ and let
$p\in \ot p.$ be a vector satisfying $p\perp P\ot q.$. 
From Proposition~\ref{Formel} and from the structure of the Fock
space $\ot F.$ (recall Eq.~(\ref{FockSpace})) we have
\begin{eqnarray*}
  p &\perp& {\rm span}\,\{a(q_1)\cdot\ldots\cdot a(q_n)\Omega\mid
             q_1,\ldots q_n \in\ot q.\,,\; n\in\N \}\;,\quad\mathrm{thus}\\
   p &\perp& \{A\Omega\mid A\in\al M.(\ot q.)\} \,.
\end{eqnarray*}
Now since $\Omega$ is cyclic for $\al M.(\ot q.)$ we must have 
$p=0$. This shows that $P\ot q.$ is dense in $\ot p.$ which is 
equivalent to $\ot p.\cap\ot q.^\perp=\{0\}$.

To show the reverse implication assume that
$\ot p.\cap\ot q.^\perp=\{0\}$. Then 
$P\ot q.$ is dense in $\ot p.$ and, consequently, the algebraic direct sum
$\mathop{\oplus}\limits_{n=0}^\infty 
 \Big(\mathop{\land}\limits^n P\ot q.\Big)$
is also dense in the antisymmetric Fock space $\ot F.$.
From Proposition~\ref{Formel} we obtain the inclusions
\[
 \mathop{\oplus}\limits_{n=0}^\infty 
 \Big(\mathop{\land}\limits^n P\ot q.\Big)
 \subset \al M.(\ot q.)\,\Omega \subset \ot F.\,,
\]
which imply that $\Omega$ is cyclic for $\al M.(\ot q.)$.
\end{proof}

Let $P$ and $Q$ be the orthoprojections corresponding to the 
subspaces $\ot p.$ and $\ot q.$ and satisfying the usual
relations w.r.t.~the antiunitary involution $\Gamma$:
\[
\Gamma P\Gamma =\EINS-P=P^\perp\quad\mr and.\quad Q\Gamma=\Gamma Q. 
\]
Proposition~\ref{Teo.3.2} above says that a necessary and sufficient condition
for doing Modular Theory with the pair $(\al M.(\ot q.),\Omega)$ is that
\begin{equation}\label{Durch1}
\ot p.\cap\ot q.=\{0\}=\ot p.\cap\ot q.^\perp\,.
\end{equation}
Using the fact that $\Gamma\ot p.=\ot p.^\perp$ we obtain, in addition,
\begin{equation}\label{Durch2}
\ot p^\perp.\cap\ot q.=\{0\}=\ot p.^\perp\cap\ot q.^\perp
           \;,\quad\mbox{where}\;\ot p.^\perp=P^\perp\ot h.\,.
\end{equation}

According to Halmos terminology (cf.~\cite{Halmos69}) if (\ref{Durch1}) and (\ref{Durch2})
hold, then the subspaces $\ot p.$ and $\ot q.$ are said to be in 
{\em generic position}. In other words the maximal subspace where 
$P$ and $Q$ commute is $\{0\}$. 
This is, in fact, a very rich mathematical situation.
For example, the following useful density statements are immediate consequences 
of the assumption that $\ot p.$ and $\ot q.$ are in generic position.
If $\ot r.\subseteq \ot q.$ (or $\ot r.\subseteq \ot q.^\perp$)
is a dense linear submanifold in $\ot q.$ (respectively in $\ot q.^\perp$),
then $P\ot r.$ is dense in $\ot p.$ and 
$P^\perp\ot r.$ is dense in $\ot p.^\perp$.
The same holds if $Q$ and $P$ are interchanged. In particular, we have that
$Q\ot p.^\perp$ is dense in $\ot q.$, $PQ\ot p.$
is dense in $\ot p.$ etc. Moreover, in the generic position situation the mapping
\begin{equation}\label{6}
P\colon \ot q. \phantom{^\perp}\longrightarrow  \ot p. 
\end{equation}
is a bounded injective linear mapping with dense image. (Similar results hold if we
interchange $\ot q.$ with $\ot q.^\perp$, $P$ with $Q$ or $Q^\perp$ etc.)

\begin{exercise}
Let $\ot q.\subset \ot h.$ be a closed $\Gamma$-invariant subspace. 
Show that the following conditions are equivalent:
\begin{itemize}
\item[(i)] $q\in \ot q.$ and $Pq=0$ implies $q=0$.
\item[(ii)] $P(\ot q.^\perp)$ is a dense submanifold of $\ot p.$.
\item[(iii)] $\ot q.\cap \ot p.=\{ 0\}$.
\end{itemize}
\end{exercise}

\begin{exercise}
Show that for any pair $P,Q$ of orthogonal projections in $\al L.(\al H.)$ one has
$\|P-Q\|\leq 1$. [Hint: Consider the operators $A:=P-Q$, $B:=\1-P-Q$ and
the equation $A^2+B^2=\1$.]
\end{exercise}

Next we give a criterion for the bicontinuity of the mapping (\ref{6}).
First, note that because $P$ is a basis projection we already have
\[
 \|PQ\|=\|QP\|=\|(\EINS-P)Q\|=\|Q(\EINS-P)\|=:\delta
\] 
and $0< \delta \leq 1$. So we can distinguish between the two
cases: $\delta < 1$ and $\delta=1$.
\begin{proposition}\label{Fall<1}
Let $P,Q$ and $\delta$ given as before. 
\begin{itemize}
\item[(i)]
If $\delta<1$, then
the mapping (\ref{6}) is bicontinuous. In particular,
the image coincides with $\ot p.$. Moreover,
the relations 
\[
 \|P-Q\|=\|(\EINS-Q)P\|=\|(\EINS-Q)(\EINS-P)\|=\delta
\]
hold.
\item[(ii)]
If $\delta=1$, then
the inverse mapping of (\ref{6}) is unbounded and densely
defined, i.e.~the image of (\ref{6}) is nontrivial 
proper dense set in $\ot p.$. 
\end{itemize}
\end{proposition}
\begin{proof}
(i) This result is a special case of Theorem~6.34 in \cite[p.~56]{bKato95}.
Note that the second alternative stated in Kato's result cannot appear in
the present situation, as a consequence of the fact that $\ot p.$
and $\ot q.$ are in generic position.

(ii) We will only show the assertion for the mapping (\ref{6}),
since one can easily adapt the following arguments to the other 
cases. Put $A:=QP^\perp Q\rest\ot q.\in\al L.(\ot q.)$, so that
$A=A^*$ and $A\geq 0$. From 
\[
 \mathrm{spr}\,A=\|A\|=\|QP^\perp P^\perp Q\|=\|P^\perp Q\|^2
 =\delta^2=1
\]
we obtain $1\in\mathrm{sp}\,A$. However, $1$ is not an eigenvalue
of $A$, because $Aq=q$, $q\in\ot q.$, implies
$\mathop{{\rm s-}\lim}\limits_{n\to \infty}\,(QP)^n q=q$
and this means $q\in\ot q.\cap \ot p.^\perp=\{0\}$. Thus
$\mathrm{ker}\,(\EINS_\ot q. -A)=\{ 0\}$ or 
$(\EINS_\ot q. -A)^{-1}$ exists and is unbounded since
$1\notin \mathrm{res}\,A$. Therefore $\vartheta:=
\mathrm{dom}\,(\EINS_\ot q. -A)^{-1}$ is a proper dense subset
in $\ot q.$ and this means $\mathrm{ima}\,(\EINS_\ot q. -A)
=\vartheta=\mathrm{ima}\,(Q -QP^\perp Q)=\mathrm{ima}\,(QPQ)$.
Finally, from the polar decomposition of $PQ$,
\[
 PQ=\mathrm{sgn}\,(PQ)\cdot(QPQ)^{\frac12}\,,
\]
we have that the partial isometry
$\mathrm{sgn}\,(PQ)$ maps $\mathrm{ima}\,(QPQ)^{\frac12}$
isometrically onto $\mathrm{ima}\,(PQ)=P\ot q.$. Thus $P\ot q.$ is
a proper dense set in $\ot p.$, i.e.~$P\colon \ot q. \to \ot p.$ is
unbounded invertible.
\end{proof}

\begin{remark}
The situation in Proposition~\ref{Fall<1}~(i) corresponds to the case where the 
index of $P$ and $Q$ is 0 (cf.~\cite[Theorem~3.3]{Avron94}).
\end{remark}

\begin{example}\label{example_proj}
As we have seen in Propositions~\ref{Fall<1},
there are two characteristic situations
when the subspaces $\ot q.$ and $\ot p.$ are in generic position.
First, when $\|PQ\|<1$. This case may be realized 
when $\ot h.$ has finite dimension.
Second, when $\|PQ\|=1$. This condition implies that the 
reference space is infinite dimensional. 
We will give here two simple examples for both situations.
\begin{itemize}
\item[(i)] {\em The case} $\|PQ\|<1$:\\
Put $\ot h.:=\C^2$ and $\Gamma(\alpha,\beta):=(\ol\beta.,\ol\alpha.)$,
$(\alpha,\beta)\in\C^2$.
The
generators of CAR$(\C^2,\Gamma)$ are simply given by
\[
 \C^2\ni (\alpha,\beta)\mapsto a(\alpha,\beta):=
  \begin{pmatrix}0 & \ol\beta. \\
        \ol \alpha. & 0
  \end{pmatrix}\;.
\]
As a basis projection we take
\[
 P:=\begin{pmatrix} 1 & 0 \\
                 0 & 0
 \end{pmatrix}\;,\quad\mr hence.\quad
\ot p.=\C \begin{pmatrix} 1\\ 0\end{pmatrix}\;,
\]
which satisfies $\Gamma P\Gamma=P^{\ms \perp.}$.
As invariant projection we choose
\[
 Q:=\frac12\begin{pmatrix} 1 & 1 \\
                           1 & 1
  \end{pmatrix}\;,\quad\mr hence.\quad
\ot q.=\C \begin{pmatrix} 1\\ 1\end{pmatrix}\;,
\]
and the $\Gamma$-invariance condition $\Gamma Q\Gamma=Q$ is trivially
satisfied.

In this example $\ot p.$ and $\ot q.$ are in generic position and it 
is straightforward to compute
\[
  \|PQ\|=\frac{1}{\sqrt2}\,.
\]

\item[(ii)] {\em The case} $\|PQ\|=1$:\\
Put $\ot h.:=L^2(\R)$ and $\Gamma f:=\ol f.$, $f\in L^2(\R)$. 
As  invariant projection define
\[
 (Q f)(x):= \chi_+(x) f(x) \;,
\]
where $\chi_+$ is the characteristic function of the 
nonnegative real numbers $\R_{\,+}=[0,\infty)$. The corresponding
$\Gamma$-invariant projection space is $\ot q.=L^2(\R_{\,+})$. 
To specify $P$ we consider first the following projection
in momentum space
\[
 (\,\wh P\, \wh f\,)(k):= \chi_+(k) \wh f(k) \;,
\]
where the Fourier transformation $F$ is defined as usual by
\[
 F(f)(k)=\wh f(k):=\frac{1}{\sqrt{2\pi}}\int_\R e^{-ikx} f(x) dx \;,
 \quad f\in L^1(\R)\cap L^2(\R)\,.
\]
Finally, the basis projection $P$ is given by
\[
 P:= F \wh P F^{-1}\,.
\]
The corresponding projection space is the Hardy space,
i.e.~$\ot p.=H^2_+(\R)$, and $P$ satisfies 
$\Gamma P\Gamma=P^\perp$.
(For a brief introduction to Hardy spaces see \cite{Baumgaertel03}).
Since by the theorem of Paley and Wiener the subspace
$H^2_+(\R)$ may be characterized in terms of holomorphic functions
on the upper half plane, it is clear that $\ot p.$ and $\ot q.$ are
in generic position. Using now the invariance of $H^2_+(\R)$ under
the regular representation $(U(a)f)(x):=f(x-a)$, $a\in\R$,
it can be shown that
\[
  \|PQ\|=1\,.
\]
\end{itemize}
\end{example}

\subsection{Modular objects for $(\al M.(\ot q.),\Omega)$:}

Let $(\al M.(\ot q.),\Omega)$ be as 
in the preceding subsection and assume that the subspace $\ot q.$ and 
$\ot p.$ are in generic position. We will compute in the
present subsection the corresponding modular objects
(recall Section~\ref{MatModularTheory}). For this purpose it is enough
to restrict the analysis to the one-particle Hilbert space $\ot p.$ of
the Fock space $\ot F.$.

Motivated by the following direct computation for the Tomita operator
\[
 S(Pq)=S\Big(a(\Gamma q)\,\Omega\Big)
  = a(\Gamma q)^*\,\Omega= a(q)\,\Omega=P\Gamma q\,,
\]
we introduce the following antilinear
mappings defined by the following graphs:
\begin{eqnarray*}
 {\rm gra}\,\beta  &:=& \Big\{(Pq,P\Gamma q)\in\ot p.\times \ot p.
                        \mid q\in\ot q.\Big\} \\
 {\rm gra}\,\alpha &:=& \Big\{(Pq^\perp,-P \Gamma q^\perp) 
                        \in\ot p.\times \ot p.
                        \mid q^\perp\in\ot q.^\perp\Big\}\,.
\end{eqnarray*}
(Note that the r.h.s.~of the preceding equations define indeed 
graphs of antilinear mappings, because the assignments
$q\to Pq$ and $q^\perp\to Pq^\perp$ are injective.)

The following result can be verified directly:
\begin{lemma}\label{Alpha}
The mappings $\alpha,\beta$ defined by the preceding graphs are
anti-linear, injective and closed with dense domains and images
$\mathrm{dom}\,\alpha=\mathrm{ima}\,\alpha=P(\ot q.^\perp)$,
$\mathrm{dom}\,\beta=\mathrm{ima}\,\beta=P\ot q.$.
Further, we have $\alpha^2=\mathrm{id}$, $\beta^2=\mathrm{id}$
on $P(\ot q.^\perp)$ resp.~$P\ot q.$ and $\alpha=\beta^*$.
\end{lemma}

\begin{remark}\label{inclusion}
The Tomita operator $S$ associated to $(\al M.(\ot q.),\Omega)$ satisfies
$S\rest\ot p.\supseteq \beta$ and $S^*\rest\ot p.\supseteq \alpha$.
Moreover, the mappings $\alpha$,$\beta$ are bicontinuous iff $\|PQ\|<1$
(recall Proposition~\ref{Fall<1}).
\end{remark}

We introduce next the notation
\[
 \Delta_\ot p.:=\beta^*\!\beta\,,
\]
since it will later turn out that $\Delta_\ot p.$ is actually 
the modular operator restricted to the one-particle Hilbert space
$\ot p.=P\ot h.$.

\begin{theorem}\label{beta*beta}
The mapping $\Delta_\ot p.\colon\ \ot p.\to \ot p.$ is a densely defined
linear positive self-adjoint operator on $\ot p.$ with graph
\[
 {\rm gra}\,\Delta_\ot p. = \Big\{(PQp,PQ^\perp p) \mid p\in\ot p.\Big\}\,.
\]
Moreover, $\Delta_\ot p.^{-1}=\beta\beta^*=\alpha^*\!\alpha$.
An expression for the modular conjugation
is given by
\[
 J(Pq)=\Delta_\ot p.^\frac12(P\Gamma q)\;.
\]
\end{theorem}
\begin{proof}
We compute first the domain of $\beta^*\!\beta$. Recalling
that $\beta^* =\alpha$ we have
\begin{eqnarray*}
\mathrm{dom}\,(\Delta_\ot p.) 
             &=&\Big\{Pq\mid q\in\ot q.\;\mathrm{and}\;
                 P\Gamma q\in\mathrm{dom}\,\alpha=P(\ot q.^\perp)\Big\}\\[2mm]
              &=&\Big\{Pq\mid q\in\ot q.\;\mathrm{and}\;
   P\Gamma q=Pq^\perp\;\mathrm{for~some}\;q^\perp\in\ot q.^\perp\Big\}\\[2mm]
             &=&\Big\{Pq\mid q\in\ot q.\;\mathrm{and}\;
                 \Gamma q\in Q(\ot p.^\perp)\Big\}\\[2mm]
             &=&\Big\{Pq\mid q\in\ot q.\;\mathrm{and}\;
                 q\in \Gamma Q(\ot p.^\perp)=Q(\Gamma \ot p.^\perp)
                 =Q\ot p.\Big\}\\[2mm]
             &=&PQ\ot p.\kern3mm= \kern3mmPQP\ot h.\,.
\end{eqnarray*}
For the third equation note that $\Gamma q-q^\perp\in\ot p.^\perp$, hence
$Q(\Gamma q)=\Gamma q\in Q(\ot p.^\perp)$. 
Since $P\Gamma Q p=-PQ^\perp\Gamma p$, 
$p\in\ot p.$ (recall $P\Gamma p=0$, $p\in\ot p.$), we have
\[
 \Delta_\ot p.(PQp)=\alpha\Big(P\Gamma Qp\Big)
 =-\alpha\Big(P Q^\perp \Gamma p\Big)=P\Gamma Q^\perp \Gamma p
 =PQ^\perp p\,,\quad p\in\ot p.\,.
\] 
Since $\ot p.$ and $\ot q.$
are in generic position the domain and image of $\Delta_\ot p.$ are dense 
in $\ot p.$. 

The last equations concerning the inverse of $\Delta_\ot p.$
follow from the preceding computation and from the fact that
$\alpha^2=\mathrm{id}$ and $\beta^2=\mathrm{id}$ on the corresponding
domains (recall Lemma~\ref{Alpha}). Finally, for the expression of the 
modular conjugation use $J=\Delta^\frac12 S$.
\end{proof}

\begin{remark}
\begin{itemize}
\item[(i)] Since $S= J\cdot\Delta^{\mt {1/2}.}$ we have from the
  preceding theorem the following inclusion of domains:
\[
 \mr dom.\Delta^{\mt {1/2}.}=\mr dom. S = P\ot q.\supset PQ\ot p.=\mr dom.\Delta\;.
\]
In this example we can characterize precisely how the domain of the
square root increases.
\item[(ii)] The present model is also useful to test many expressions
  that appear in general computation done in Modular Theory. For
  example, for certain calculations one needs to work with the dense
  set $\al D.:=\mr dom.\Delta^{\mt {1/2}.}\cap\mr dom.\Delta^{-\mt
    {1/2}.}$.  In the present example involving the CAR-algebra it is
  straightforward to verify that $\al D.=PQ\ot p.^\perp$, which is in
  fact dense in $\ot p.$ since the corresponding two subspaces are in
  generic position.
\end{itemize}
\end{remark}
\begin{exercise}
Show that, in general, for a positive self-adjoint operator $T$ in a 
Hilbert space one has the inclusion: 
\[
 \mr dom. T^\frac12\supseteq \mr dom. T\,.
\]
\end{exercise}

We conclude mentioning the behavior of the modular objects 
with respect to the direct sums that appear in the Fock space
$\ot F.$ (cf.~Eq.~(FockSpace)). For a complete proof see \cite{Lledo02a}. 
Let $(\ot h.,\Gamma)$, $P$ and
$Q$ be as before and denote by $S=J\Delta^\frac12$ the polar decomposition 
of the Tomita operator for the pair $(\al M.(\ot q.),\Omega)$. 

Note that the different modular objects leave the 
$n$-particle submanifolds $\mathop{\land}\limits^n(P\ot q.)$ invariant. 
(This fact is well known in the context of CCR-algebras 
\cite{pLeyland78}, where one can use the so-called 
exponential vectors which are specially well-adapted to the 
Weyl operators.)

\begin{proposition}\label{NTeilchen}
Let $q_1,\ldots,q_n\in \ot q.$ and $q^\perp_1,\ldots,
q^\perp_n\in \ot q.^\perp$. Then the following equations hold
\begin{eqnarray*}
S(P q_1\land \ldots \land Pq_n)
    &=& P \Gamma q_n\land \ldots\land P\Gamma q_1\kern2mm=\kern2mm
        S(P q_n)\land \ldots \land S(Pq_1 ) \\
S^*(P q_1^\perp\land \ldots \land Pq_n^\perp)
    &=& S^*(P q_n^\perp)\land \ldots \land S^*(Pq_1^\perp)
\end{eqnarray*}
Moreover, 
\[
  \mathrm{span}\,\Big\{ a(q_1)\cdot\ldots\cdot a(q_n)\Omega\mid
                 q_1,\ldots,q_n\in\ot q.\;,\;n\in\N\cup\{0\}\Big\}
\]
is a core for the Tomita operator $S$. 
\end{proposition}

The following result together with Theorem~\ref{beta*beta} gives a complete
picture of the modular objects in the context of CAR-algebras.
\begin{theorem}\label{ModRest}
Let $(\al M.(\ot q.),\Omega)$ be as in the preceding 
subsection and assume that $\ot q.$ and 
$\ot p.$ are in generic position.
Let $S=J\Delta^\frac12$ be the polar decomposition of the Tomita operator.
The modular operator $\Delta =S^*S$ and the modular conjugation $J$ 
can be restricted to the respective $n$-particle subspaces. 
In particular we have:
\begin{itemize}
\item[(i)] Modular operator: We have $\Delta\rest\ot p.=\Delta_\ot p.$, 
  where $\Delta_\ot p.=\beta^*\!\beta$, and
  $\mathrm{dom}\,\Delta\rest P_n\ot F.=\mathop{\land}\limits^n
  \mathrm{dom}\,\Delta_\ot p.$. Moreover, the action on 
  the $n$-particle vector is given by
 \begin{equation}\label{MOn}
  \Delta (p_1\land\ldots\land p_n)=(\Delta_\ot p.\,p_1)
  \land\ldots\land(\Delta_\ot p.\,p_n)\,,\quad 
  p_1,\ldots,p_n\in\mathrm{dom}\,\Delta_\ot p.=PQ\ot p.\,.
 \end{equation}
\item[(ii)] Modular conjugation: Its action on 
  the $n$-particle vector is given by
\begin{equation}\label{MCn}
J (p_1\land\ldots\land p_n)=(Jp_n)
  \land\ldots\land (Jp_1)\,,\quad p_1,\ldots,p_n\in\ot p.\,.
\end{equation}
\end{itemize}
\end{theorem}

\begin{remark}
The CAR-algebra is typically
used to model Fermi systems in quantum physics, while bosonic
systems are described in terms of the CCR-algebra.
A formula for the modular operator was given for the (bosonic) free 
scalar field in \cite{Figliolini89}. In this paper the reference
space is specified in terms of the Cauchy data 
of the Klein-Gordon operator and the formula
for the modular operator on the one-particle Hilbert space reads
\[
  \delta=\frac{B+1}{B-1}\,,
\]
where the operator $B$ is defined in terms of two other 
densely defined closed operators $A_{\mt {\pm 1}.}$ and these
are again defined using suitable idempotents $P_{\mt {\pm 1}.}$
(see \cite[p.~425]{Figliolini89} for details).

The simplicity of the formulas obtained in the context of the 
self-dual CAR-algebra (see e.g.~Theorem~\ref{beta*beta})
suggest that also for the bosonic models the self-dual approach 
to the CCR-algebra may be better adapted to problems concerning 
Modular Theory. In fact, in this case one can also characterize
the Fock representations in terms of basis projections 
(cf.~\cite{Araki71/72a,Araki71/72b,bPetz90}). Therefore, it seems
likely that similar simple formulas as the ones presented in 
this chapter also hold in the context of the CCR-algebra.
\end{remark}

\subsection{Modular objects for double cones in Fermi models}
We mention finally that the formulas established previously
also apply to the localized algebras that appear
in the context of Fermi free nets (see e.g.~\cite{Lledo01,Lledo04}
and references therein). For more details on local quantum theories
see Section~2 in \cite{pGuido08}.
Let $\al O.\subset\R^4$ be a double cone in Minkowski space
and denote by $\overline{\ot q.(\al O.)}$ the closure of the subspaces 
$\ot q.(\al O.)$ of the reference Hilbert space $(\ot h.,\Gamma)$. 
The subspaces $\ot q.(\al O.)$ are defined 
in terms of the embeddings that characterize
the free nets (essentially Fourier transformation of $C^\infty$ functions
with compact support restricted to the positive mass shell/light cone).
It is easily shown that $\Gamma\ot q.(\al O.)=\ot q.(\al O.)$,
hence $\Gamma\overline{\ot q.(\al O.)}=\overline{\ot q.(\al O.)}$.
Moreover the localized C*-algebras are again CAR-algebras:
\begin{eqnarray*}
 \al A.(\al O.)&:=&\mathrm{C}^*\left(\{a(\varphi)\mid
                 \varphi\in\ot q.(\al O.) \}\right)
               \;=\;\mathrm{CAR}(\ot q.(\al O.),\Gamma\restriction\ot q.(\al O.)) \\
               &=&\mathrm{CAR}(\overline{\ot q.(\al O.)},
                \Gamma\restriction\overline{\ot q.(\al O.)})
                \subset\mathrm{CAR}(\ot h.,\Gamma)\,,
\end{eqnarray*}
where for the last equation we have used Proposition~\ref{StabRef}.
For the canonical basis projection $P$ given in the context of 
Fermi free nets (see e.g.~\cite[p.~1157]{Lledo01}) and 
for double cones $\al O.$ one has
\[
 \ot p.\cap\overline{\ot q.(\al O.)}=
 \ot p.\cap\overline{\ot q.(\al O.)}^\perp=\{0\}\;,\quad\mathrm{where}\;
 \ot p.=P\ot h.\,,
\]
(see also Section~2 in \cite{pGuido08} for similar relations in the 
case of Bose fields).
Therefore $\ot p.$ and $\overline{\ot q.(\al O.)}$ are in generic 
position and we can apply the results and formulas of the
present section to the von Neumann algebras
in a Fock representation specified by $P$ and localized in double cones $\al O.$:
\[
\al M.(\al O.):=\left\{ a(\varphi)\mid \varphi\in\overline{\ot q.(\al O.)}\right\}''
                \,. 
\]
In particular from Proposition~\ref{NTeilchen} and
Theorem~\ref{ModRest} we have that the modular operator $\Delta$ 
and the modular conjugation $J$ are already characterized by their
action on the one-particle Hilbert space. Finally, using
Theorem~\ref{beta*beta} we have:

\begin{theorem}\label{O}
Let $\al O.\subset\R^4$ be a double cone in Minkowski space.
Denote by $Q_\al O.$ the orthoprojection onto $\overline{\ot q.(\al O.)}$
and by $P$ the canonical basis projection given in the context of Fermi
free nets. Then the following formulas hold for the modular operator and
modular conjugation on the one-particle Hilbert space $\ot p.$:
\begin{eqnarray*}
  \mathrm{gra}\,\Delta_\ot p. 
        &=&  \Big\{(PQ_\al O.(p),PQ_\al O.^\perp (p)) \mid p\in\ot p.\Big\}\,.\\
  J(Pq) &=& \Delta_\ot p.^{\frac12}(P\Gamma q)\;,
            \quad q\in\overline{\ot q.(\al O.)}\,.
\end{eqnarray*}
\end{theorem}

\section{Some classical applications of Modular Theory}

In this section we mention briefly some classical applications of Modular
Theory to the theory of von Neumann algebras. For further applications
in quantum field theory see \cite{pGuido08}.

\subsection{The commutant of tensor products}\label{subsec:commutants}

Recall that in Example~\ref{ExFinBg} we have used the modular
conjugation $J$ to show that
\[
 (\al L.(\al H.)\otimes\C\1_{\al H.'})'=\C\1_{\al H.}\otimes\al L.(\al H.')\,.
\]
This is just a special case of the following general result:
\begin{theorem}
Let $\al M.$, $\al N.$ be von Neumann algebras. Then
\[
 (\al M.\otimes\al N.)'=\al M.'\otimes\al N.'\;.
\]
\end{theorem}
In this generality the statement of the theorem remained open for a
long time and the first proof used Modular theory.

\subsection{Structure of type~III factors}\label{subsec:typeIII}

The technically more tractable cases mentioned previously in
Examples~\ref{ex:LeftReg} and \ref{ExFinBg}, as well as in
Section~\ref{sec:ex-crossed}, have in common that the corresponding
von Neumann algebra $\al M.$ is finite or, equivalently, that the
vector state associated with the cyclic and separating vector is a
trace. In order to treat infinite algebras one has to consider more
general states or even weights\footnote{Let $\al M.$ be a von Neumann
  algebra and denote by $\al M._+$ its positive elements. A weight
  $\varphi$ on $\al M.$ is a map $\varphi\colon\al M._+\to [0,\infty]$
  which is additive and positively homogeneous.}.

To get deeper into the structure of type~III factors it is necessary
to consider Modular Theory in the more general context defined by
Hilbert Algebras and focus on the crucial information contained in the
modular automorphism group. Recall that the action of the modular
automorphism group is nontrivial if the von Neumann algebra is
infinite.  In this more general context one can also associate to any
faithful, normal, semifinite weight $\varphi$ on a von Neumann algebra
$\al M.$ modular objects $(\Delta^\varphi,J^\varphi)$ (for details see
\cite{bSunder87,bKadisonII,bTakesakiII}). Connes analyzed in
\cite{Connes73} (see also \cite{ConnesIn76} for a review or
\cite{bTakesaki83,bSunder87}) the dependence of the modular
automorphism group $\sigma^\varphi_t$ on the weight $\varphi$. In
\cite{Connes73} the author established the following fundamental
theorem:
\begin{theorem}
Let $\varphi,\psi$ be faithful, normal, semifinite weights on the von Neumann algebra
$\al M.$. Then there is a $\sigma$-weakly continuous one-parameter family 
$\{U_t\}_{t\in\R}$ of unitaries in $\al M.$ satisfying the cocycle condition
\[
 U_{t+s}=U_t\sigma_t^\varphi(U_s)\;,\quad t,s\in\R\;,
\]
and such that 
\[
 \sigma_t^\psi(M)=U_t\;\sigma_t^\varphi(M)\;U_t^*\;,\quad M\in\al M.\;.
\]
\end{theorem}

This theorem shows that the modular automorphism group is essentially (up to 
unitary equivalence) independent of the weight and any von Neumann algebra 
carries a natural action of $\R$ given by the modular automorphism group.
Connes also introduced in his seminal paper \cite{Connes73} the following
algebraic invariant associated to the von Neumann algebra $\al M.$ and that
uses the spectrum of the modular operators:
\[
  \Gamma(\al M.):=\bigcap \Big\{\mathrm{sp}(\Delta^\varphi)\mid
                        \varphi~\mathrm{faithful,~normal,~semifinite~weight~on~}
                        \al M. \Big\}\,.
\]
This invariant is crucial for the finer classification of type~III factors.
In fact, it turns out that $\Gamma(\al M.)$ is a closed multiplicative semigroup
of $\R_+$ and, therefore, there are only the following possibilities:
\[
\Gamma(\al M.)=\left\{
\begin{array}{l}
[0,\infty) \\[2mm] 
\{\lambda^n\mid n\in\Z\}\cup\{0\}\quad\mathrm{for~some~}\,\lambda\in(0,1) \\[2mm]
\{0,1\}\;.
\end{array}
\right. 
\]
In either case one says that the type~III factor $\al M.$ is of type~III$_0$,
of type~III$_\lambda$, $\lambda\in(0,1)$, or of type~III$_1$, respectively.
This classification is unique if the factor is hyperfinite, i.e.~if it is 
generated by an increasing sequence of finite-dimensional *-subalgebras.
\subsection{The KMS condition}\label{subsec:KMS}

As was seen in the preceding subsection the modular automorphism group
$\{\sigma_t\}_{t\in\R}$ plays a fundamental role in the classification
of type~III factors.  There is also a characteristic and very useful
analytic relation between $\{\sigma_t\}_{t\in\R}$ and the
corresponding state $\varphi$. For simplicity, we will formulate it in
the case where $\varphi$ is a faithful normal {\em state.}
\begin{definition}
  A one-parameter automorphism group $\{\sigma_t\}_{t\in\R}$ satisfies
  the modular condition relative to the state $\varphi$ if invariance holds,
  i.e.~$\varphi\circ\sigma_t=\varphi$, $t\in\R$, and if for each
  $M,N\in\al M.$ there is a complex-valued function $F$ satisfying the
  following two conditions:
\begin{itemize}
\item[(i)] $F$ bounded and continuous on the horizontal strip
  $\{z\in\C\mid 0\leq\mathrm{Im}(z)\leq 1\}$ and analytic on the
  interior of that strip.
\item[(ii)] $F$ satisfies the following boundary condition:
\[
 F(t)=\varphi(\sigma_t(M)N) \quad\mathrm{and}\quad
 F(t+i)=\varphi(N\sigma_t(M))\;,\quad t\in\R.
\]
\end{itemize}
\end{definition}

To each state $\varphi$ as above there corresponds uniquely a one-parameter automorphism
group $\sigma_t$ satisfying the modular condition. In view of the preceding 
definition, the modular automorphism group gives a measure of the extent to which
the state fails to be tracial (see also Remark~\ref{rem-mod}~(iii)).

\begin{remark}
The modular condition mentioned before is known in quantum statistical mechanics
as the KMS (Kubo-Martin-Schwinger) boundary condition 
(in this case at inverse temperature $\beta=-1$). In this context
$\{\sigma_t\}_{t\in\R}$ describes the time evolution of the system
and KMS condition was proposed as a criterion for equilibrium
(see \cite{bHaag92} for further details).
\end{remark}

\section{Appendix: Crossed products and the CAR-algebra}

In this appendix we collect some material used in the examples
presented in Sections~\ref{sec:ex-crossed} and \ref{sec:CAR} 
and that would have interrupted the flow of the article.

\subsection{Crossed products}\label{subsec:crossed}

Recall the group measure space construction presented
in Section~\ref{sec:ex-crossed}. In particular let
$(\Omega,\Sigma,\mu)$ be a $\sigma$-finite measure space
and denote by $\al M.:=L^\infty(\Omega,\Sigma,\mu)$ the corresponding
maximal abelian von Neumann algebra in $\al L.(\al H.)$, where
$\al H.:=L^2(\Omega,\Sigma,\mu)$. Denote by 
$\alpha\colon \Gamma \to\mathrm{Aut}\cM$ the action of the 
discrete group $\Gamma$ on $\al M.$ and by 
$\al N.=\al M.\otimes_\alpha\Gamma$ the corresponding crossed 
product acting on $\al K.:=\ell_2(\Gamma)\otimes\al H.$.

\begin{proposition}
With the preceding notation we have:
\begin{itemize}
\item[(i)] The image $\pi(\al M.)$ of $\al M.$ in the crossed product
 $\al N.$ is maximal abelian iff the action $\alpha$ is free (cf.~Eq.~(\ref{abelianM})).
\item[(ii)] Assume that the action $\alpha$ is free. Then $\alpha$ is ergodic iff
 $\al N.$ is a factor.
\end{itemize}
\end{proposition}

The next theorem shows that
all types of von Neumann factors mentioned before may be 
realized explicitly within the group measure space construction
described previously. It implies, in particular, Theorem~\ref{teo:finite}.

\begin{theorem}\label{ex-factors}
Let $(\Omega,\Sigma,\mu)$ and 
$\al M.:=L^\infty(\Omega,\Sigma,\mu)$ be as in the preceding proposition.
Assume that there is a free and ergodic action 
$\alpha\colon\Gamma\to\mathrm{Aut}\al M.$ of a discrete group
$\Gamma$ on $\al M.$. For the types of the factor
$\al N.=\al M.\otimes_\alpha\Gamma$ we have the following criteria:
\begin{itemize}
\item[(i)] Suppose that there is $\Gamma$-invariant $\sigma$-finite positive
measure $\nu$ which is equivalent to $\mu$ (in the sense of mutual absolute
continuity). Then
\begin{itemize}
    \item[$\bullet$] $\al N.$ is of type~I iff the measure space
      $(\Omega,\Sigma,\mu)$ contains atoms.
    \item[$\bullet$] $\al N.$ is of type~II iff the measure space
      $(\Omega,\Sigma,\mu)$ contains no atoms.
    \item[$\bullet$] $\al N.$ is finite iff $\nu$ is a finite measure.
\end{itemize}
\item[(ii)] The factor $\al N.$ is of type~III iff there does not exist a $\sigma$-finite
      positive measure $\nu$ which is equivalent to $\mu$ and $\Gamma$-invariant.
\end{itemize}
\end{theorem}

More concrete and explicit examples (including type~III$_0$, type~III$_\lambda$
and type~III$_1$ factors)
can be found in \cite[{\S}4.3]{bSunder87},
\cite[{\S}V.7]{bTakesakiI},
\cite[{\S}XIII.1]{bTakesakiIII} or \cite[{\S}8.6]{bKadisonII}.

\subsection{The self-dual CAR-algebra}\label{CARalgebras}
In this subsection we recall some standard results on the self-dual
CAR-algebra which is needed in the example in Section~\ref{sec:CAR}.
We will define and state the main properties of the C*-algebra that
is associated to the canonical anticommutation relations and its 
irreducible representations. General references for the
present section are \cite{Araki70/71,ArakiIn87}.  

\begin{theorem}
\label{Teo.2.2.9}
Let $\ot h.$ be a complex Hilbert space with scalar product
$\langle\cdot,\cdot\rangle$ and anti-unitary involution $\Gamma$,
i.e.~$\langle\Gamma f,\Gamma h\rangle=
\langle h,f \rangle$, for all $f,h\in \ot h.$. Then
$\car{\ot h.,\Gamma}$ denotes
the algebraically unique C*-algebra generated by
$\EINS$ and $a(\varphi)$, $\varphi\in \ot h.$, such that the following
relations hold:
\begin{itemize}
\item[{\rm (i)}] The mapping 
$\ot h.\ni\varphi\mapsto a(\varphi)$ is antilinear.
\item[{\rm (ii)}] $a(\varphi)^*=a(\Gamma\varphi)\;$, $\varphi\in \ot h.\,$.
\item[{\rm (iii)}] $a(\varphi_1) a(\varphi_2)^*+a(\varphi_2)^*
  a(\varphi_1)=\langle\varphi_1,\varphi_2\rangle \;\EINS\;$,
  $\varphi_1,\varphi_2\in \ot h.\,$.
\end{itemize}
\end{theorem}

$\ot h.$ is called the {\em reference space} of 
$\car{\ot h.,\Gamma}$. This space is a `parameter' space 
labeling the generators of the algebra.
The name `self-dual' comes from property 
(ii) above, where the algebra involution $^*$ is described
in terms of the antilinear mapping $\Gamma$ of the 
reference space $\ot h.$. For a finite-dimensional
example see Example~\ref{example_proj}~(i).
The preceding uniqueness result implies the following 
statement concerning the automorphisms of the CAR-algebra.

\begin{theorem}
\label{Teo.2.2.10}
Let $U$ be a unitary of the reference space $\ot h.$ that satisfies
$U\,\Gamma=\Gamma\, U$. Any such $U$ generates an automorphism 
$\alpha_{\KIn{$U$}}$ of $\car{\ot h.,\Gamma}$ $($called the {\em Bogoljubov
automorphism} associated to the {\em Bogoljubov unitarity} $U)$ uniquely 
determined by the equation
\[
  \alpha_{\KIn{$U$}}(a(\varphi)):= a(U\varphi)\,,\quad\varphi\in \ot h.\,.
\]
\end{theorem}

\begin{definition}
\label{Def.2.2.11}
An orthoprojection $P$ on the reference space $\ot h.$ is a {\em basis
projection} if it satisfies the equation $P+\Gamma P\Gamma=\EINS$.
\end{definition}

\begin{theorem}
\label{Teo.2.2.12}
Any basis projection $P$ generates a unique state
$\omega_{\KIn{$P$}}$ on $\car{\ot h.,\Gamma}$ by means of the relation
\[
 \omega_{\KIn{$P$}}\!\LR a(\varphi)^*a(\varphi) \RR=0,\quad\Wort{if} 
 \quad P\varphi=0\,.
\]
$\omega_{\KIn{$P$}}$ is a pure state and is called the {\em Fock state}
corresponding to the basis projection $P$.
\end{theorem}

An explicit representation $\pi$ of the CAR-algebra associated to a 
basis projection $P$ is realized on the
the antisymmetric Fock space
\begin{equation}\label{FockSpace}
 \ot F.:= \mathop{\oplus}\limits_{n=0}^\infty 
          \Big(\mathop{\land}\limits^n P\ot h.\Big)\,.
\end{equation}
At this point we introduce the usual annihilation and creation
operators on $\ot F.$.
\begin{eqnarray*}
c(p)\,\Omega    &:=& 0\,,\\
c(p)\,(p_1 \land \ldots \land p_n) 
                &:=& \sum\limits_{r=1}^{n}\, (-1)^{r-1}\, 
              \langle p,p_r\rangle_{\got h} \;p_1 
              \land\ldots\land\widehat p_r \land \ldots \land p_n\,,\\[2mm]
c(p)^*\,\Omega  &:=& p\,,\\
c(p)^*\,(p_1 \land \ldots \land p_n) 
                &:=& p \land p_1 \land \ldots \land p_n\,,
\end{eqnarray*}
where $\Omega$ is the Fock vacuum in the subspace corresponding to
$n=0$ in the definition (\ref{FockSpace}) and $p,p_1,\ldots, p_n\in
P\ot h.$. The symbol 
$\widehat p_r$ means that the vector $p_r$ is omitted
in the (antisymmetric) wedge product $\land$. 
Finally, the Fock representation $\pi$
is defined by
\[
 \pi(a(f)):= c(P\Gamma f)^*+c(Pf)\,,\quad f\in \ot h.\,.
\]

In the rest of this section we assume that a basis projection $P$ 
is given and when no confusion arises we will also simply write $a(f)$ 
instead of $\pi(a(f))$. We will later need an explicit expression for
$a(f_n)\cdot\ldots\cdot a(f_1)\,\Omega$. Let $n,k,p$ be natural numbers
with $2p+k=n$ and define the following subset of the symmetric group
$\ot S._n$:
\[
\begin{array}{l} \ot S._{n,\,p}:=\left\{\left( {\begin{array}{cccccccc}
       n & n-1     &\cdots & n-2p+2  & n-2p+1 & k   & \cdots & 1 \\
\alpha_1 &\beta_1 &\cdots &\alpha_p &\beta_p & j_1 & \ldots & j_k   
               \end{array}}\right)\in \ot S._n\; \right. \\[5mm] 
 \;\;\;\;\;\;\;\;\;\;\;\;\;\;\;\;\; \left. \Big|\;\;
 \alpha_1>\ldots>\alpha_p\,,\; \alpha_l>\beta_l\,,\;l=1,\ldots, p 
 \quad\mbox{and}\quad n\geq j_1>j_2>\ldots j_k\geq 1
 \right\}\,.\end{array}
\]
Note that $\ot S._{n,\,p}$ contains $\left(\kern-1.5mm\begin{array}{c}
n \\n-2p \end{array}\kern-1.5mm\right) {\displaystyle \frac{(2p)!}{
 p!\,2^p}}\;$ elements.

\begin{proposition}\label{Formel}
For $f_1,\ldots,f_n\in \ot h.$ the equation
\[
 \Big(a(f_n)\cdot\ldots\cdot a(f_1)\Big)\,\Omega = 
  \sum\limits_{\mbox{\tiny $\begin{array}{c}\pi\in\ot S._{n,\,p}
               \\[1mm] 0\leq 2p\leq n\end{array}$}}
  \!\!\!\!({\rm sgn}\,\pi)\;
  \prod\limits_{l=1}^p \;\langle Pf_{\alpha_l}\,,\,P\Gamma f_{\beta_l}\rangle
  \, P\Gamma f_{j_1} \land \ldots \land P\Gamma f_{j_k} 
\]
holds, where the indices $\alpha_l, \beta_l,j_1,\ldots,j_k$ are given
in the definition of $\ot S._{n,\,p}$ and where for $n=2p$ in the preceding
sum one replaces the wedge product by the vacuum $\Omega$.
\end{proposition}

Finally, we state the following proposition that shows the 
stability of the CAR-algebra w.r.t.~the operation of
taking the closure of the reference space. 
\begin{proposition}\label{StabRef}
Let $\ot h._0$ be a complex pre-Hilbert space and $\Gamma_0$ an 
antilinear involution on it. Denote by $(\ot h.,\Gamma)$ the corresponding
closures. Then
\[
 \car{\ot h._0,\Gamma_0}:=\mr C.^*\Big(a(\varphi)\mid \varphi\in \ot h._0\Big)
                   =\car{\ot h.,\Gamma}\,,
\]
where ${\rm C}^*(\cdot)$ denotes the C*-closure of the argument.
\end{proposition}

\providecommand{\bysame}{\leavevmode\hbox to3em{\hrulefill}\thinspace}

\end{document}